\documentstyle[12pt]{article}
\newcommand{\be}{\begin{equation}}
\newcommand{\ee}{\end{equation}}
\newcommand{\ba}{\begin{array}}
\newcommand{\ea}{\end{array}}
\newcommand{\ben}{\begin{enumerate}}
\newcommand{\een}{\end{enumerate}}

\newcommand{\ov}{\overline}
\newcommand{\tr}{{\rm tr}\, }
\newcommand{\ad}{{\rm ad}\, }
\newcommand{\id}{{\rm id}\, }

\newcommand{\w}{{\!}\wedge{\!}}

\newcommand{\im}{\mbox{Im}\, }

\font \msb=msbm10 scaled \magstep1

\newcommand{\bR}{\mbox{\msb R} }
\newcommand{\bC}{\mbox{\msb C} }
\font \msbm=msbm7 scaled \magstep0
\newcommand{\bRm}{\mbox{\msbm R} }
\newcommand{\bCm}{\mbox{\msbm C} }

\font \msbl=msbm9 scaled \magstep1
\newcommand{\bCl}{\mbox{\msbl C} }

\newcommand{\br}{\beta }

\newcommand{\dr}{\delta }
\newcommand{\er}{\varepsilon }
\newcommand{\lr}{\lambda }

\newcommand{\cg}{\ov{g}}

\newcommand{\sw}{s^{\wedge}}
\newcommand{\gd}{g^{\dagger}}

\newcommand{\dc}{\partial}
\newcommand{\bc}{\ov{\partial}}

\font \eul=eufm10 scaled \magstep2

\newcommand{\gotG}{\mbox{\eul g}}
\newcommand{\gotH}{\mbox{\eul h}}

\begin{document}

\title{\bf Free motion on the Poisson $SU(N)$ group}
\author{{\bf S. Zakrzewski}  \\
\small{Department of Mathematical Methods in Physics,
University of Warsaw} \\ \small{Ho\.{z}a 74, 00-682 Warsaw, Poland} }

\date{}
\maketitle

\begin{abstract}
$SL(N,\bCl )$ is the phase space of the Poisson $SU(N)$. We
calculate explicitly the symplectic structure of $SL(N,\bCl )$,
define an analogue of the Hamiltonian of the free motion on
$SU(N)$ and solve the corresponding equations of motion.
Velocity is related to the momentum by a non-linear Legendre
transformation.
\end{abstract}

\section{Introduction}

The theory of Poisson groups \cite{D:ham,D,S-T-S,Lu-We,Lu} (and
their phase spaces \cite{CDW,Ka,qcp,srni})
allows us to consider deformations of known mechanical models.
Several low-dimensional examples, related to Poisson symmetry,
have been already investigated
\cite{poi,poican,k-part,zakop,standr,sphere}. All these models
certainly have their quantum-mechanical counterpart, the
underlying Poisson group being replaceable by the corresponding
quantum group. It is natural to study first the Poisson case as
technically simpler. We obtain interesting classical systems,
and, at the same time, we get some idea about the corresponding
quantum systems.

In this paper we calculate explicit form of Poisson brackets on
the phase space of the Poisson $SU(N)$ group, i.e. on $SL(N,\bC
)$ (as a real manifold). We also consider a natural candidate
for the Hamiltonian of the free motion. It turns out that the
projections of the phase trajectories onto $SU(N)$ are `big
circles' (shifted one-parameter subgroups), as in the usual
case. The (constant) velocity is however a non-linear function
of the momentum, so we have an example of a deformed Legendre
transformation.

The case of $SU(2)$ was presented in \cite{zakop}, in a direct
(tedious) way --- without referring to the compact $r$-matrix
notation. The above mentioned deformed character of the Legendre
transformation in this case was shown in \cite{sphere} to be the
reason why the free dynamics reduced to the homogeneous space
(Poisson sphere) yields really a deformation of usual free
trajectories on the sphere.

The paper is organized as follows. In Section~2 we clarify when
the Drinfeld double of a Lie bialgebra $(\gotG ,\dr)$  coincides
with the complexification of $\gotG $ (recall \cite{Lu-We} that
this is the case of $\gotG = su(N)$), and we obtain a useful
formula for the Drinfeld's canonical $r$-matrix $r_D$ on the
double. In Section~3 we calculate $r_D$ for $\gotG = su(N)$ in
terms of matrix units.  This allows effectively to write down
the Poisson brackets of matrix elements of $SL(N,\bC )$. In
Section~4 we introduce the free Hamiltonian, which is one of the
most natural functions on $SL(N,\bC )$. We solve the equations
of motion and analyze the bijectivity property of the `Legendre
transformation'.

\section{Drinfeld double and complexification}

Let $\gotG$ be a real Lie algebra and let $b(\cdot ,\cdot )$ be
a non-degenerate invariant symmetric bilinear form on $\gotG$.
On the complexification $\gotG ^{\bCm}$ we have then the
non-degenerate invariant symmetric bilinear form $B:=\im b^{\bCm}$,
with respect to which both $\gotG$ and $i\gotG$ are isotropic.
We are thus almost in the situation of a Manin triple: all
properties are satisfied except that $i\gotG $ is not
a Lie subalgebra (unless $\gotG $ is abelian).

Of course, any isotropic Lie subalgebra $\gotH$ in
$\gotG ^{\bCm}$, which is complementary to $\gotG$, yields a
Manin triple and the corresponding Lie bialgebra structure on
$\gotG$. The question now arises, which Lie bialgebra structures
on $\gotG$ are obtained in this way.

A subspace $\gotH$ of $\gotG ^{\bCm}$ is complementary to
$\gotG$ if it is of the form
\be\label{iR}
 \gotH = \{ ix+Rx : x\in \gotG \},
\ee
where $R\colon\gotG\to\gotG$ is a linear map. Such a subspace is
isotropic (with respect to $B$) if and only if $R$ is
skew-symmetric with respect to $b$:
\be\label{skew}
 b(Rx,y)=-b(x,Ry),\qquad x,y\in\gotG .
\ee
This subspace is a subalgebra if and only if
\be\label{mod}
[Rx,Ry]-R([Rx,y]+[x,Ry])=[x,y],\qquad x,y\in \gotG .
\ee
Let $s\in\gotG\otimes \gotG$ denote the inverse of $b$. We shall
use the same letter for $s$ considered as a linear map from
$\gotG ^*$ to $\gotG$. The composition  $r:=Rs$ is then a
linear map from $\gotG ^*$ to $\gotG $, which we can also
identify with an element of $\gotG\otimes\gotG$ (an element of
$\gotG\otimes\gotG$ defines a linear map from $\gotG ^*$ to $\gotG$
by the contraction in the first argument). In terms of this
$r\in\gotG \otimes \gotG $,
condition (\ref{skew}) means that $r$ is anti-symmetric, and
condition (\ref{mod}) is equivalent to
\be\label{rmod}
[[r,r]]=[[s,s]],
\ee
where $[[w,w]]$ for $w\in \gotG\otimes\gotG$ denotes the Drinfeld's
bracket
$$
[[w,w]]:=[w_{12},w_{13}]+[w_{12},w_{23}]+[w_{13},w_{23}].
$$
In terminology of \cite{standr}, it means that $r+is$ is an {\em
imaginary} quasi-triangular classical $r$-matrix:
$$
[[r+is,r+is]]=0.
$$
It is easy to show that the Lie bialgebra structure on $\gotG$
defined by $r$ (by taking the coboundary of $r$), coincides with
the one defined by the Manin triple $(\gotG ^{\bCm}; \gotG ,\gotH)$:
$$
B([ix+Rx,iy+Ry],z)= b([Rx,y]+[x,Ry],z)
$$
$$
=b(y,\ad _z Rx)-b(x,\ad _z Ry)=
(b\otimes b)(x\otimes y,(\id\otimes \ad _z)r)+
(b\otimes b)(x\otimes y,(\ad_z\otimes \id )r)
$$
$$
=(B\otimes B)((ix+Rx)\otimes (iy+Ry),\ad _z r).
$$
We have thus established a one-to-one correspondence between
Manin triples realized in $\gotG ^{\bCm}$ (with the scalar
product $B$) and imaginary quasi-triangular Lie
bialgebra structures on $\gotG$ (with the
symmetric part of the $r$-matrix equal $s=b^{-1}$).
The Drinfeld double of this Lie bialgebra is $\gotG ^{\bCm}$
with the quasi-triangular structure given by the canonical element
\be\label{wD}
 w_D = e_k\otimes f^k\in \gotG ^{\bCm}\otimes \gotG ^{\bCm}
\ee
(summation convention), where $e_k$ is a basis of $\gotG$ and
$f^j$ is the dual (w.r.t. $B$) basis in $\gotH$. One
can easily see that $f^j=r(e^j)+is(e^j)$, where $e^j$ is the
dual basis in $\gotG ^*$.

The skew-symmetric part $r_D$ of $w_D$ is given by
\be\label{rD}
r_D= \frac12 e_j\wedge [r(e^j)+is(e^j)]=r + \frac12 e_j\wedge
(is^{jk}e_k).
\ee
Note that both $w_D$ and $r_D$ are elements of the real tensor
product $V\otimes _{\bRm}V$, where $V:=\gotG ^{\bCm}$ is treated
as a real vector space. They may be, however treated as (real)
elements of the complexification
$$
(V\otimes _{\bRm}V)^{\bCm}=V^{\bCm}\otimes _{\bCm}V^{\bCm},
$$
which is much more convenient. In order to distinguish the
imaginary unit $i$ arising in the complexification of $V$ from
the imaginary unit arising in the complexification of $\gotG$,
we denote the latter by $J\colon \gotG\to \gotG$. Recall, that
any $v\in V$ may be represented as the sum
$$
v = v^+ + v^-,\qquad v^{\pm }:=\frac12 (v\pm \frac{1}{i}Jv),
$$
and we have
$$
(Jv)^{\pm } = \pm iv^{\pm }\qquad (\mbox{we have also}\;\;\;
[v_1^{\pm},v_2^{\pm}]=\pm [v_1,v_2]^{\pm }).
$$
In a fixed basis, it is also convenient to set
$$
e_j=\dc _j + \bc _j, \qquad \mbox{where}\;\;
\dc _j:= e_j^+,\; \bc _j:=e_j^-.
$$
In particular, the last term in (\ref{rD}) may be written as follows:
$$
\frac12 s^{jk}(\dc_j + \bc _j)\wedge (i\dc_j - i\bc_k)=
is^{jk}\bc_j \w \dc_k.
$$
This term will be denoted by $\sw $. Note that
\be\label{s+-}
\sw = \w (\id \otimes J)s,
\ee
where $\w (a\otimes b):=a\w b$ is the antisymmetrization.

\section{Drinfeld and Heisenberg double of Poisson $SU(N)$}

Let $b$ denote the invariant scalar product on $\gotG:=su(N)$ given by
\be\label{bXY}
b(X,Y):= - \frac{1}{\er } \tr XY
\ee
($\er$ is a parameter).
Let $\gotH =sb(N)$ be the Lie subalgebra in $\gotG
^{\bCm}=sl(N,\bC )$ consisting of complex uppertriangular
matrices with real diagonal elements (and trace zero).
It is easy to see that $\gotH$ is complementary to $\gotG$ and
isotropic with respect to $B=\im b^{\bCm}$, hence $(\gotG
^{\bCm};\gotG,\gotH)$ is a Manin triple. It corresponds to the
standard Poisson $SU(N)$ \cite{Lu-We}. Our aim is to calculate
$r_D= r + \sw $ given by (\ref{rD}). We introduce first the
typical elements of $su(N)$ defined in terms of usual matrix
units $e_j{^k}=e_j\otimes e^k$:
$$F_j{^k} := e_j{^k} - e_k {^j},\qquad G_j{^k} := i(e_j{^k} + e_k
{^j}),\qquad H_{jk} :=i( e_j{^j} - e_k{^k})
$$
so that
\be\label{basis}
F_j{^k}, \; G_j{^k}\;\; (j<k), \qquad H_{j,j+1}\;\; (1\leq j\leq
N-1)
\ee
is a basis of $su(N)$.

\vspace{2mm}

\noindent
{\bf Lemma 1.} \ {\em We have}
\be\label{rsun}
r=\frac{\er}{2} \sum_{j<k} F_j{^k}\w G_j{^k},
\ee
\be\label{ssun}
s =  \frac{\er}{2} \sum_{j<k} (F_j{^k}\otimes F_j{^k}+
G_j{^k}\otimes G_j{^k} +\frac{2}{N}H_{jk}\otimes H_{jk} ).
\ee

\vspace{1mm}

\noindent
{\em Proof:} \ It is easy to calculate first $R$ defined in
(\ref{iR}). Since $ix+Rx\in sb(N)$ for $x\in su(N)$, it is
easy to calculate the lowertriangular part of $Rx$ (it is the
corresponding part of $-ix$) and the diagonal part of $Rx$ (it is
the diagonal part of $-ix$ plus something real, hence zero). We
obtain $RH_{jk}=0$, $RF_j{^k}=G_j{^k}$, $RG_j{^k}=-F_j{^k}$.
Since $F_j{^k}$, $G_j{^k}$ $(j<k)$ form an orthogonal set with
$$ b(F_j{^k},F_j{^k})=\frac{2}{\er }=b(G_j{^k},G_j{^k}),$$
and they are orthogonal to all $H_{jk}$,
it is easy to check that contraction of the $r$ given in
(\ref{rsun}) with the basis elements (using $b$) coincides with
the action of $R$ on these elements.

In order to prove (\ref{ssun}), note first that only the last
term needs to be explained (due to the orthogonality). The
scalar product on the Cartan subalgebra spanned by $H_{jk}$ is
naturally the restriction of the scalar product defined on the
space spanned by $H_j:=ie_j{^j}$ (the Cartan for $u(N)$) by the
same formula:
$$ b(H_j,H_k)=-\frac{1}{\er}\tr H_jH_k =\frac{1}{\er} \dr _{jk}.$$
In order to invert $b$ on the subspace, it is sufficient to
invert it on the bigger space, which is easy:
\be\label{HjHj}
 \er \sum_j H_j\otimes H_j,
\ee
and project it orthogonally on the subspace. Since
$$ H_j = \frac1{N} \sum_k H_k +\frac1{N}\sum_k (H_j-H_k)$$
is the orthogonal decomposition, we just have to replace $H_j$
in (\ref{HjHj}) by $\frac1{N}\sum_k (H_j-H_k)$. This indeed
gives (\ref{ssun}).

\hfill $\Box $

\vspace{1mm}

Now we can express $r$, $\sw $ and finally $r_D$ in terms
of `holomorphic' and `anti-holomorphic' vectors $\dc _j{^k}=
(e_j{^k})^+$, $\bc _j{^k}=(e_j{^k})^-$. This will enable us to
calculate the Poisson brackets of basic coordinate functions
(and their complex conjugates) on $SL(N,\bC )$.

A straightforward insertion of $e_j{^k}=\dc _j{^k}+\bc _j{^k}$,
$Je_j{^k}=i\dc _j{^k}-i\bc _j{^k}$
into (\ref{rsun}) and (\ref{ssun}) together with (\ref{s+-}) gives
$$ r = r^{(2,0)}+r^{(0,2)}+r^{(1,1)},$$
where
\be
r^{(2,0)}=\ov{r^{(0,2)}}= i\er \sum_{j<k} \dc_j {^k}\w \dc_k{^j},\qquad
r^{(1,1)}= i\er \sum_{j<k} (\bc_j {^k}\w \dc_j{^k} -\bc_k {^j}\w \dc_k{^j}),
\ee
and
\be
\sw =-i\er (\frac1{N}\ov{I}\w I -\sum_{j,k}\bc _j{^k}\w \dc _j{^k}),
\ee
where $ I:=\sum_k\dc _k{^k}$.

The $r$-matrix $r_D= r+\sw $ on $sl(N,\bC )$ defines two
Poisson bivector fields on $SL(N,\bC )$:
$$ \pi _{\pm} (g) = r_D (g\otimes g) \pm (g\otimes g)r_D,\qquad
g\in SL(N,\bC )$$
{\em Drinfeld double} of the Poisson $SU(N)$ is the Poisson group
$(SL(N,\bC ),\pi _-)$. The {\em Heisenberg double} of the
Poisson $SU(N)$ is the Poisson manifold $(SL(N,\bC ),\pi _+)$.
It plays the role of the phase space (cotangent bundle) of the
Poisson $SU(N)$. The bivector field $\pi _+$ is known to be
non-degenerate (\cite{Lu,qcp}), because $SL(N,\bC )$ globally
decomposes (by the Iwasawa decomposition) onto $G=SU(N)$ and
$G^*:=SB(N)$, i.e. every element $g\in SL(N,\bC )$ is a product
of the form
$$  g = u\br , \qquad u\in G,\; \br \in G^* $$
with uniquely defined $u,\br $. Here $SB(N)$ is the connected
subgroup of $SL(N,\bC )$, corresponding to the Lie algebra
$\gotH =sb(N)$ (i.e. the Poisson {\em dual} of the Poisson
$SU(N)$).

Using the compact notation
$$ \{ g_1,g_2\} ^{ab}_{cd} = \{ g^a_c,g^b_d\},\qquad
(g_1g_2)^{ab}_{cd}=(g\otimes g)^{ab}_{cd}=g^a_cg^b_d ,$$
we can now write the Poisson brackets of matrix elements of $g$
for $\pi _{\pm }$ as follows
\be\label{sunpb}
 \{ g_1,g_2\}_{\pm} = \rho g_1g_2\pm g_1g_2\rho ,\qquad
\{ \cg_1,g_2\}_{\pm} = w' \cg_1g_2\pm \cg_1g_2w' ,
\ee
where $\rho :=r_D^{(2,0)} = r^{(2,0)}$ is the purely holomorphic
part of $r_D$ and
$$
w':= -i\er (\frac1{N}\ov{I}\otimes I -\sum_k \bc_k{^k}\otimes
\dc _k{^k} -2\sum_{j< k} \bc _j{^k}\otimes \dc _j{^k}) .
$$
is the antiholomorphic-holomorphic part (without
antisymmetrization) of $r_D$, i.e. $r_D^{(1,1)}=
r^{(1,1)}+\sw =w'-w'_{21}$.

The Poisson structure on $SL(N,\bC )$ viewed as the Drinfeld double of
the Poisson $SU(N)$ is therefore described by the brackets
\begin{eqnarray*}
\{ g^j_l,g^j_m\}_-  & = &  i\er g^j_lg^j_m \;\;\;\;\;\; (l<m) \\
\{ g^j_l,g^k_l\}_-  & = &  i\er g^j_lg^k_l \;\;\;\;\;\;\; (j<k) \\
\{ g^j_l,g^k_m\}_-  & = &  2i\er g^j_mg^k_l \;\;\;\; (l<m,j<k) \\
\{ g^j_l,g^k_m\}_-  & = &  0 \;\;\;\;\;\;\;\;\;\;\;\;\;\;\; (l>m,j<k)
\end{eqnarray*}
as far as the holomorphic variables are concerned (quite
standard), and 
\be\label{mix}
\{ \cg^j_l,g^k_m\}_-  = i\er [\dr ^{jk}(\cg^j_lg^j_m +
2\sum_{a>j}\cg^a_lg^a_m ) - \dr _{lm} (\cg^j_lg^k_l
+2\sum_{a<l}\cg^j_ag^k_a ) ]
\ee
for the mixed case (this is nothing but the Poisson version of
commutation relations for the real quantum group $SL(N,\bC )$,
cf. \cite{Po}, formulae (3.77)--(3.80)).

The Heisenberg double Poisson structure on $SL(N,\bC )$ is given
by
\begin{eqnarray*}
\{ g^j_l,g^j_m\}_+  & = & - i\er g^j_lg^j_m \;\;\;\;\;\; (l<m) \\
\{ g^j_l,g^k_l\}_+  & = &  i\er g^j_lg^k_l \;\;\;\;\;\;\; (j<k) \\
\{ g^j_l,g^k_m\}_+  & = & 0 \;\;\;\;\;\;\;\;\;\;\;\;\;\;\; (l<m,j<k) \\
\{ g^j_l,g^k_m\}_+  & = &   2i\er g^j_mg^k_l \;\;\;\;  (l>m,j<k)
\end{eqnarray*}
\be\label{mix+}
\{ \cg^j_l,g^k_m\}_+  = i\er [-\frac2{N} \cg^j_lg^k_m +
\dr ^{jk}(\cg^j_lg^j_m + 2\sum_{a>j}\cg^a_lg^a_m ) 
+ \dr _{lm} (\cg^j_lg^k_l + 2\sum_{a<l}\cg^j_ag^k_a ) ] .
\ee

It is sometimes convenient to replace the complex conjugate
variable $\cg$ by $\gd$ --- the hermitian conjugate of $g$.
In this case the second equality in (\ref{sunpb}) is replaced by
\be\label{gdg0}
\{ \gd _1,g_2\}_{\pm} = \gd _1 [(\tau \otimes \id)w']g_2\pm
g_2 [(\tau \otimes \id)w']\gd _1,
\ee
where $\tau $ denotes the transposition. If we set
$$ w := -(\tau \otimes \id)w'= i\er (\frac1{N} I\otimes I-\sum_k
\dc_k{^k}\otimes \dc _k{^k} -2\sum_{j< k} \dc _k{^j}\otimes \dc _j{^k}),
$$
we can write these brackets as follows:
\be\label{gdg}
\{ \gd _1,g_2\}_{\pm} = -\gd _1 wg_2\mp g_2 w\gd _1 .
\ee
Note that $\rho = \frac12 (w - w_{21})$ is the antisymmetric
part of $w$,
\be\label{w-r}
w -\rho = \frac12 (w+w_{21}) = i\er (\frac1{N}I\otimes I - P),
\ee
where $P$ is the permutation, and $iw$ is the infinitesimal part
of the $R$-matrix  for the $A_N$-series (cf. \cite{FRT,Po}),
$${\cal R}=I\otimes I +iw+\ldots $$
(when $q=1+\er +\ldots $), hence  $w$ satisfies the classical
Yang-Baxter equation $[[w,w]]=0$.

\section{Free motion on Poisson $SU(N)$}

The Poisson structure on $SL(N,\bC )$ viewed as the Heisenberg
double of Poisson $SU(N)$ (analogue of the cotangent bundle) is
given by
\be
\{ g_1,g_2\} = \rho g_1g_2+ g_1g_2\rho ,\qquad
\{ \gd _1,g_2\} = -\gd _1 wg_2- g_2 w\gd _1
\ee
(we have dropped the subscript `+', for simplicity).

In the non-deformed case of the cotangent bundle $T^*G$ to
$G=SU(N)$, the free motion is governed by the Hamiltonian
$H\colon T^*G\to \bR$ proportional to the square of
the momentum, given by a biinvariant metric on $G$. In other words,
there is a distinguished quadratic function on the dual $\gotG
^*$ of the Lie algebra $\gotG$ of $G$ (defined by the Killing
form), and $H$ is just the pullback of this function to $T^*G$
(from the left or from the right --- it does not matter since
the quadratic function is invariant under coadjoint action of
$G$; also: it is a Casimir for the Poisson structure on $\gotG^*$).

When $T^* SU(N)$ is replaced by $SL(N,\bC )$, it is still
easy to find a Hamiltonian with similar properties, namely
\be
H (g) := \frac12 \tr \gd g .
\ee
Note that it depends only on the `momenta' $\br \in G^*$:
$$
H(g)=H(u\br ) = \frac12 \tr (u\br )^{\dagger} u\br =
\frac12 \tr \br ^{\dagger}\br = H(\br ).
$$
It does not matter which way we decompose $g$:
$$
H(g)=H(\br ' u')=\frac12 \tr (\br ' u')^{\dagger} \br ' u' =
\frac12 \tr {\br '}^{\dagger }\br ' = H(\br ' ),
$$
which means that $H$ as a function on $G^*$ is invariant with
respect to the dressing action:
$$ H(\br )=H(u\br ) = H( {^u\br}\cdot u^{\br })= H({^u\br })$$
(notation of \cite{qcp,sphere}). It means that $H$ is a Casimir
on $G^*$.

We shall now examine the equations of motion. We have
$$ \dot{g} = \{ H,g\} = \frac12 \tr _1 \{ \gd _1g_1,g_2\} 
=\frac12 \tr _1 (\{ \gd _1,g_2\}g_1 + \gd _1\{ g_1,g_2\} )
$$
$$
= \frac12 \tr_1 (-\gd _1 wg_2g_1 -g_2w\gd _1g_1 + \gd _1\rho
g_1g_2 +\gd_1 g_1g_2\rho )=-\frac12 \tr_1 [\gd _1(w-\rho )g_1g_2
+ \gd _1g_1g_2 (w-\rho )].
$$
Using (\ref{w-r}) and the identity
$$ \tr _1 \gd _1g_1g_2P = g\gd g = \tr _1 g_1\gd _1 Pg_2, $$
we obtain
\be\label{eqmot}
 \dot{g} = i\er [g\gd g - \frac1{N} (\tr \gd g) g].
\ee
Substituting here $g=u\br $, we get
$$ \dot{u}\br + u\dot{\br} = i\er [u\br \br
^{\dagger}u^{\dagger} u\br - \frac1{N} (\tr \br ^{\dagger}\br)
u\br ],$$
or,
$$ u^{-1}\dot{u} + \dot{\br}\br ^{-1} = i\er [\br\br ^\dagger
-\frac1{N} (\tr \br\br ^{\dagger})].$$
Since the right hand side belongs to $\gotG = su(N)$, we have
$\dot{\br} =0$, which was in fact also clear before, because $H$
is a Casimir on $G^*$. Therefore we are left
with the condition of constant velocity
\be
u^{-1}\dot{u}  = F(\br ) :=i\er [\br\br ^\dagger 
- \frac1{N} (\tr \br\br ^{\dagger})].
\ee
It follows that as far as configurations are concerned, the
motion looks exactly as the non-deformed one: the particle moves
on the `big circles' (shifted 1-parameter subgroups) with
constant velocity. The difference consists in the momentum
variables, which have a non-linear nature. The function $F$
above tells how to compute the velocity when the momentum is
given. It plays the role of the inverse Legendre transformation.

A general notion of the Legendre transformation in the case of
phase spaces of Poisson manifold is investigated in \cite{tLt}.
Here we shall show only two properties of the map $F\colon
SB(N)\to su(N)$:
\ben
\item it intertwines the dressing action with the adjoint action:
$$  F({^u}\br ) = uF(\br ) u^{-1},$$
\item it is bijective.
\een
The first property follows from the fact that if $u\br =\br
'u'$, then $\br '\br'^{\dagger} = u\br\br ^{\dagger}u^{\dagger}$.
To prove the second, we first show that the map
$$ SB(N)\ni \br \mapsto \psi (\br )= \br\br ^{\dagger}\in P:= \{ p :
p>0,\; \det p =1\} $$
is a bijection. Define a map $\phi : P\to SB(N)$ by
$$ \phi (p) := \br,\qquad \mbox{where $\br$ is such that} \;\;
p^{\frac12} = \br u\;\; (\mbox{Iwasawa}).$$
We have $\psi \circ \phi =\id$, since $\br\br ^{\dagger}=(\br
u)(\br u)^{\dagger} = p^{\frac12}p^{\frac12} =p$.
But $\phi $ is also surjective, since, given $\br\in SB(N)$, it
is sufficient to consider its polar decomposition $\br =
p_0u_0$ and notice that $\phi (p_0^2)=\br$.

It remains to prove that the map
$$ P\ni p\mapsto h= p-\frac1{N} \tr p \in \, i\cdot su(N)$$
is a bijection. We first show that $h$ determines $p$. Choose an
orthonormal basis in which $h$ is diagonal, then $p$ is also
diagonal in that basis. Let $p_i$ and $h_i$ be the corresponding
eigenvalues, then
$$ \lr _i = p_i -\left\langle p\right\rangle, \qquad
\mbox{where} \;\; \left\langle p\right\rangle := \frac1{N}
\sum_j p_j.$$
If $\lr _i$ come from some $p$, then
$$ \lr _i + \left\langle p\right\rangle >0,\qquad
(\lr _1 + \left\langle p\right\rangle )\cdot \ldots\cdot (\lr _N
+ \left\langle p\right\rangle ) =1.
$$
Since the function
$$
[\max_j (-\lr _j),\infty [ \; \ni t\mapsto f(t):=
(\lr _1 +t)\cdot\ldots\cdot (\lr _N +t)\in [0,\infty [
$$
is a (monotonic) bijection, there is exactly one $t_0$ such that
$f(t_0)=1$, hence $\left\langle p\right\rangle =t_0$ and this
completely determines $p$ by $p_i=\lr _i
+\left\langle p\right\rangle $. It is easy to see that $p_i =
\lr _i + t_0$, where $f(t_0) =1$, defines some $p\in P$ for every
$h$ (because then $p_i>0$ and $p_1\cdot\ldots\cdot p_n=1$).

Finally, we remark that in the limit $\er\to 0$, the model
becomes the undeformed one:
$$ \br \sim I +\er \xi,\qquad \xi\in  sb(N)\equiv su(N)^*
,\qquad F(\br )\sim i(\xi +\xi ^{\dagger} ),
$$
and
$$
\frac{H(\br ) -\frac1{N}}{\er ^2} \sim \frac12 \tr \xi\xi
^{\dagger }.
$$
In view of the existence of a Poisson isomorphism between
$SB(N)$ and $su(N)^*$ \cite{GW}, it would be interesting to
find, how the function $H$ on $SB(N)$ is expressed as a function
on $su(N)^*$.

\end{document}